\documentclass{aptpub}

\usepackage{float}
\usepackage{amsmath}
\usepackage{amsfonts}
\usepackage{amssymb}
\usepackage[T1]{fontenc}
\usepackage{geometry}
\usepackage{graphicx}
\usepackage{hyperref}
\usepackage{icomma}
\usepackage[latin1]{inputenc}
\usepackage{latexsym}
\usepackage[numbers]{natbib}
\usepackage{textcomp}
\usepackage[all]{xy}
\usepackage{cancel}
\usepackage{color}
\usepackage{dsfont}
\usepackage{subfig}
\usepackage{caption}

\authornames{MART\'INEZ S., SAN MART\'IN J., VILLEMONAIS D.}
\shorttitle{QSD and fast return from infinity}

\DeclareMathSymbol{\minus}{\mathord}{operators}{"2D}

\def\be{\begin{eqnarray}}
\def\ee{\end{eqnarray}}
\def\ben{\begin{eqnarray*}}
\def\een{\end{eqnarray*}}
\def\bei{\begin{itemize}}
\def\eei{\end{itemize}}
\def\me{\medskip \noindent}
\def\bi{\bigskip \noindent}
\def\E{\mathbb{E}}
\def\P{\mathbb{P}}

\def\1{\mathbf{1}}
\def\N{\mathbb{N}}

\def\d{\partial}

\begin{document}

\title{Existence and uniqueness of a quasi-stationary distribution for Markov processes with fast return from infinity}

\authorone[Universidad de Chile]{Servet Mart\'inez}
\addressone{Departemento Ingenier\'ia, Matem\'atica and Centro de Modelamiento Matem\'atico, Universidad de Chile, UMI 2807, CNRS -- Universidad de Chile}
\authortwo[Universidad de Chile]{Jaime San Mart\'in}
\addresstwo{Departemento Ingenier\'ia, Matem\'atica and Centro de Modelamiento Matem\'atico, Universidad de Chile, UMI 2807, CNRS -- Universidad de Chile}
\authorthree[University of Lorraine]{Denis Villemonais}
\addressthree{Institut \'Elie Cartan de Lorraine, University of Lorraine, France}

\begin{abstract}
We study the long time behaviour of a Markov process evolving in $\N$ and conditioned not to hit $0$. Assuming that the process comes back quickly
from infinity, we prove that the process admits a unique \textit{quasi-stationary distribution} (in particular, the distribution of the conditioned process admits a limit when time goes to infinity). Moreover, we prove that the distribution of the process converges exponentially fast in total variation norm to its quasi-stationary distribution and we provide an explicit rate of convergence.

As a first application of our result, we bring a new insight on the speed of convergence to the quasi-stationary  distribution for birth and death processes: we prove that these processes converge exponentially fast to a quasi-stationary distribution if and only if they have a unique quasi-stationary distribution. Also, considering the lack of results on quasi-stationary distributions for non-irreducible processes on countable spaces, we show, as a second application of our result, the existence and uniqueness of a quasi-stationary distribution for a class of possibly non-irreducible processes.
\end{abstract}

\keywords{process with absorption; quasi-stationary distributions;
Yaglom limit; mixing property; birth and death processes}

\ams{37A25; 60B10; 60F99}{60J80}

\section{Introduction}
Let $X$ be a stable continuous time Markov process evolving in
$\N=\{0,1,2,...\}$ such that $0$ is an absorbing point, so
$X_t=0\ \forall t\geq T_0$, and absorption occurs almost surely, that is
for all $x\in\N$, $\P_x(T_0<+\infty)=1$, where $T_0=\inf\{s\geq
0,\ X_s=0\}$.
In this paper, we
provide a sufficient condition for the existence and uniqueness of
a quasi-stationary distribution for $X$ and for the conditional distribution of $X$ to converge exponentially fast to it. 

\me A \textit{quasi-stationary distribution} (QSD) for
$X$ is a probability measure $\rho$ on $\N^*=\{1,2,3,...\}$ such
that, for all $t\geq 0$,
\begin{equation*}
  \rho(\cdot)=\P_\rho\left(X_t\in\cdot|t<T_0\right).
\end{equation*}
Thus a QSD is stationary for the process conditioned not to be absorbed.
The notion of QSD has always been closely
related to the study of the long-time behavior of a process conditioned not to
be absorbed. Indeed, it is well known (see for instance~\cite{VereJones1969}, \cite{Meleard2011}) that a probability measure
$\rho$ is a QSD if and only if it is a \textit{quasi-limiting distribution}
(QLD), which means that there exists a probability measure $\mu$ on $\N^*$ such that
\begin{equation}
\label{eq:QLD}
 \rho(\cdot)= \lim_{t\rightarrow\infty}\P_\mu\left(X_t\in\cdot | t<T_0\right).
\end{equation}
%
Existence and uniqueness of QSDs and QLDs have
been extensively studied in the past decades. They have originally
been investigated by Yaglom~\cite{Yaglom1947}, which stated their existence for sub-critical Galton-Watson processes.
  In their seminal work~\cite{Darroch1967}, Darroch
and Seneta proved that irreducible finite state space processes admit a unique QSD.  In our case of a
process $X$ evolving in a countable state space, the question is more
intricate since the existence or uniqueness of a QSD is not always true.
In 1995,
Ferrari, Kesten, Mart\'inez and Picco~\cite{Ferrari1995} proved a
necessary and sufficient condition for the existence of a
quasi-stationary distribution for $X$ under the assumption that it is
irreducible and that the process doesn't come back from infinity in
finite time. More precisely, the authors proved that if $\N^*$ is an
irreducible class for the process $X$ and if
$\lim_{x\rightarrow+\infty} \P_x(T_0<t)=0$ for any
$t>0$, then the existence of a QSD for $X$ is equivalent to $
\E_x\left(e^{\lambda T_0}\right)<+\infty$ for some constants $x\in\N^*$ and $\lambda>0$. 
The much--studied birth and death processes are of particular interest, since explicit sufficient and necessary conditions have been proved by
van Doorn~\cite{vanDoorn1991} characterising the three possible cases:  there is no QSD, a unique QSD or an infinite continuum of QSDs. (For more informations on QSDs/QLDs, we refer the reader to the recent surveys~\cite{Meleard2011} and~\cite{vanDoorn2011}.) In this paper, we give a sufficient criterion for the existence and uniqueness of a quasi-stationary distribution for countable state space processes. In the particular case of birth and death processes, we shall see that the criterion is in fact equivalent to the existence and uniqueness of a QSD.

\me While the existence of a QSD is interesting in itself, it is only the first step towards the understanding of a conditioned process long time behaviour. Indeed, it is of first practical importance to determine the initial distributions $\mu$ for which the convergence~\eqref{eq:QLD} holds and, as stressed out in~\cite{Meleard2011}, to determine the speed of convergence to the QSD. 
 In the present paper, our aim is twofold since we give a criterion ensuring  the existence and uniqueness of a QSD and we prove that the conditional distribution of the process converges exponentially fast in total variation norm to a unique quasi-stationary distribution. Moreover, we provide an explicit speed of convergence, independent of the initial distribution of the process.  
More precisely, we prove that there exist a unique QSD $\rho$ and a constant $\gamma\in]0,1[$ (for which we provide an explicit expression) such that
\begin{equation*}
\|\P_\mu(X_t\in\cdot|t<T_0)-\rho(\cdot) \|_{TV} \leq 2(1-\gamma)^{[t]},\ \forall \mu\in{\cal M}_1(\N^*),\ \forall t\geq 0,
\end{equation*}
where $\|\cdot\|_{TV}$ denotes the total variation norm for signed measures, $[t]$ is the integer part of $t$ and ${\cal M}_1(\N^*)$ refers to the set of probability measures on $\N^*$. As we shall see, our proof uses a purely probabilistic approach, allowing us to answer the long standing question of the speed of convergence of a birth and death process to its unique quasi-stationary distribution, which was out of reach of the spectral theory tools historically used to handle this case. 

\me The existence and uniqueness criterion is based on the three following hypotheses, where the positive constants $c_1, c_2, c_3$ and $c_4$ will appear in the expression of $\gamma$. Our first assumption H1 states that there exists a subset of $\N^*$ where the probability of extinctions at any time $t$ are balanced. 
\begin{description}
\item[Hypothesis H1] There exists a finite subset of $K\subset\N^*$ and a constant $c_1>0$ such that, for all $t\geq 0$,
      \begin{equation*}
        \frac{\inf_{x\in K}\P_{x}(t<T_0)}{\sup_{x\in K}\P_x(t<T_0)}\geq c_1.
      \end{equation*}
\end{description}
\textit{Remark.}
One easily check that when the process is irreducible, that is when $\P_x(X_t=y>0)$ for all $x,y\in \N^*$, this property is fulfilled for any finite subset $K\subset \N^*$. Note also that the smallest the subset $K$ is, the weakest the requirement on the constant $c_1$ is.

\me Let $K$ satisfying H1. Our second assumption is that the process comes back quickly from any point to $K\cup\{0\}$ and, starting from some particular point in $K$, it has a relatively high probability to be in $K$ afterwards. We denote by $T_K=\inf\{n\geq 0,\;X_n\in K\}$ the hitting time of $K$. 
\begin{description}
\item[Hypothesis H2] There exist some constants $\lambda_0>0, c_2>0, c_3>0$ and a point $x_0\in K$ such that, for all $t\geq 0$,
      $$\sup_{x\in \N^*} \E_x(e^{\lambda_0 T_K\wedge T_0})\leq c_2\text{ and }\P_{x_0}\left(X_t\in K\right)\geq c_3 e^{-\lambda_0 t}.$$
\end{description}
\textit{Remark.} Usually, there exists an interval of values of $\lambda_0$ acceptable here, as it will clearly appear in the birth and death case (see the proof of Theorem~\ref{theorem:BD} below).
 Note also that the largest the subset $K$ is, the weakest the requirements on the constants $\lambda_0, c_2,c_3$ are.

\me Our last assumption is that the conditioned process comes back in time $1$ to a point $x_0\in K$ with a minimal probability.
\begin{description}
\item[Hypothesis H3] There exists a constant $c_4>0$ and a point $x_0\in K\subset E$, such that
       $$\inf_{x\in\N^*} \P_x(X_1=x_0\ |\ T_0>1)\geq c_4.$$
\end{description}
\textit{Remark.} If the rate of absorption is uniformly bounded over $\N^*$, then $\inf_{x\in\N^*}\P_x(T_0>1)>0$ and thus assumption H3 is equivalent to the existence of $x_0$ and $c_4$ such that $\inf_{x\in\N^*} \P_x(X_1=x_0)\geq c_4$. This is closely related to the existence of a \textit{small set}, following the terminology of Down, Meyn and Tweedye~\cite{DownMeynTweedie1995} for processes without absorption, where $T_0=+\infty$ happens $\P_x$-almost surely.

\me We are now able to state our main theorem, which is proved in Section~\ref{section:main-result}. As an application, we also provide a corollary on birth and death processes and show a generalization of the
recent results of Ferrari and Mari\`c~\cite{Ferrari2007}.

\begin{theorem}
\label{theorem:main}
  If Hypotheses~H1, H2 and~H3 are fulfilled, then
there exists a unique QSD $\rho$ for $X$. Moreover, for any
  probability measure $\mu$ on $\N^*$, we have
  \begin{align}
  \label{equation:mixing:in-theorem:main}
    \|\P_\mu(X_t\in\cdot|t<T_0)-\rho \|_{TV} \leq 2 \left(1-\frac{c_1c_3c_4}{2c_2}\right)^{[t]},\ \forall t\geq 0.
  \end{align}
\end{theorem}

\me\textit{Remark.} Inequality~\eqref{equation:mixing:in-theorem:main} implies that $\rho$ is a QLD for $X$ and any initial distribution,
  which means that, for any probability measure $\mu$ on $\N^*$,
  \begin{equation*}
    \lim_{t\rightarrow\infty} \P_\mu(X_t\in \cdot | t< T_0)=\rho(\cdot).
  \end{equation*}

\me\textit{Remark.} Our approach is based on a strong mixing
property inspired
by~\cite{DelMoralVillemonais2011}. In particular, we prove that
\begin{equation*}
\|\P_\mu(X_t\in\cdot|t<T_0)-\P_\nu(X_t\in \cdot |t<T_0) \|_{TV} \leq 2\left(1-\frac{c_1c_3c_4}{2c_2}\right)^{[t]},\ \forall \mu,\nu\in{\cal M}_1(\N^*),\ \forall t\geq 0.
\end{equation*}

\me\textit{Remark.} Notice that $c_1,c_2,c_3,c_4$ can be chosen in a way that they satisfy $\frac{c_1c_3c_4}{2c_2}<1$. Nevertheless,
as a consequence of the proof of Theorem 1 below these constants always satisfy $\frac{c_1c_3c_4}{2c_2}\le 1$ (see the argument after equation (\ref{eq:ZZZ}).

\me We present two applications of our result. In  Section~\ref{section:application1}, we develop the case of birth and death processes and prove that such a process admits a unique quasi-stationary distribution $\rho$ if and only if Theorem~\ref{theorem:main} holds, that is if its conditional distribution converges exponentially fast to $\rho$, uniformly in its initial distribution. Note that this result  provides a very new insight on the quasi-limiting behaviour of birth and death processes. Moreover, its proof reveals that our criterion is optimal for birth and death process: such a process satisfies Hypotheses H1, H2 and H3 if and only if it admits a unique QSD.

 \me In our second
application, developed in Section~\ref{section:application2}, we show that the sufficient condition for existence and
uniqueness of a QSD proved in~\cite{Ferrari2007} can be considerably
relaxed. While the practical implications of this application is nowadays less manifest than the previous one, it is of much theoretical interest. Indeed, it demonstrates that our result applies to reducible Markov processes on a countable state space, which is an exciting area under development where most of the existing results on QSDs do not apply.

\section{Proof of Theorem~\ref{theorem:main}}
\label{section:main-result}

The proof of Theorem~\ref{theorem:main} is divided into three
parts. In a first step, we show that, for all $t\geq 0$,
\begin{equation}
  \label{equation:c4}
  \frac{\P_{x_0}\left(t<T_0\right)}{\sup_{x\in \N^*} \P_x\left(t<T_0\right)}\geq \frac{c_1 c_3}{2c_2}.
\end{equation}
Secondly, using the techniques developed in Del Moral and
Villemonais~\cite{DelMoralVillemonais2011}, we prove
inequality~\eqref{equation:mixing:in-theorem:main} for all $t\geq
0$. In a third step, we conlude the proof by showing that~\eqref{equation:mixing:in-theorem:main} implies the existence and
uniqueness of a QSD.

\bi Step 1: Let us show that~\eqref{equation:c4} holds.  
For all $x\in E$, we have
  \begin{align*}
    \P_x(t<T_0) &= \E_x\left( \1_{t < T_K \wedge T_0} \right) +\E_x\left(\1_{T_K \leq t <T_0}\right).
  \end{align*}

\me On the one hand, we deduce from Hypothesis~H2 that, for all $t\geq 0$,
\begin{align*}
  \E_x\left( \1_{t < T_K \wedge T_0} \right) &\leq e^{-\lambda_0 t}
  \E_x\left(e^{\lambda_0 T_{K}\wedge T_0}\right)
  \leq  \frac{\P_{x_0}(X_t\in K)}{c_3} c_2\\
  &\leq \frac{\P_{x_0}(X_t\in \N^*)}{c_3} c_2=\frac{\P_{x_0}(t<T_0)}{c_3} c_2
\end{align*}

\me On the other hand, the Markov property yields to
\begin{align*}
  \E_x\left(\1_{T_K \leq t <T_0}\right)&= \E_x\left( \1_{T_{K}\leq t\leq
    T_0} \P_{X_{T_{K}}}(t-T_{K}\leq T_0) \right)\\
  &=  \E_x\left(\1_{T_{K}\leq t\leq T_{0}}
    e^{\lambda_0 T_{K}\wedge T_0} e^{-\lambda_0 T_{K}\wedge T_0}
    \P_{X_{T_{K}}}(t-T_{K} < T_0)\right)\\
   &\leq \E_x\left(e^{\lambda_0 T_{K}\wedge T_0}\right) \sup_{y\in
      K} \sup_{s\in[0,t]} e^{-\lambda_0 s}\P_y(t-s < T_0)\\
   &\leq c_2 \sup_{y\in
      K} \sup_{s\in[0,t]} e^{-\lambda_0 s}\P_y(t-s < T_0),
\end{align*}
by Hypothesis~H2.  Now, by
Hypotheses~H1 and~H2, we have for all
$s\in[0,t]$ and any $y\in K$,
  \begin{align*}
    e^{-\lambda_0 s}\times \P_{y}(t-s < T_0) 
    &\leq  \frac{\P_{x_0}(X_{s}\in K)}{c_3}\times \frac{\inf_{z\in K} \P_z(t-s < T_0)}{c_1}\\
    &\leq \frac{\P_{x_0}(t < T_0)}{c_1 c_3},
  \end{align*}
  where we used the Markov property. We deduce that
  \begin{equation*}
   \E_x\left(\1_{T_K \leq t <T_0}\right)\leq \frac{c_2}{c_1 c_3} \P_{x_0}(t<T_0).
  \end{equation*}

\me Finally, we have
\begin{equation*}
  \P_x(t<T_0)\leq \left( \frac{c_2}{c_3} + \frac{c_2}{c_1 c_3} \right)\P_{x_0}(t<T_0)
\end{equation*}
  which implies~\eqref{equation:c4}, since $c_1$ is necessarily smaller than $1$.

\bi Step 2: Let us define, for all $0 \leq s \leq
 t\leq T$ the linear operator $R_{s,t}^T$ by
\begin{align*}
  R_{s,t}^T f(x)&= \E_x(f(X_{t-s})\mid T-s< T_0)\\
                &= \E(f(X_{t})\mid X_s=x,\ T< T_0),
\end{align*}
by the Markov property. We begin by proving that, for any $T>0$, the
family of operators $(R_{s,t}^T)_{0\leq s\leq t\leq T}$ is a Markov semi-group. We
have, for all $0\leq u\leq s\leq t\leq T$,
 \begin{equation*}
   R_{u,s}^T (R_{s,t}^T f)(x) = \E_x(\E_{X_{s-u}}(f(X_{t-s})\ |\ T-s< T_0) \ |\ T-u< T_0).
 \end{equation*}
For any measurable function $g$, the Markov property implies that
 \begin{align*}
   \E_x\left(g(X_{s-u})\1_{T-u< T_0}\right)
   &=\E_x\left(g(X_{s-u})\P_{X_{s-u}}(T-u -(s-u) < T_0)\right)\\
   &=\E_x\left(g(X_{s-u})\P_{X_{s-u}}(T-s < T_0)\right).
 \end{align*}
Applying this equality to $g:y\mapsto \E_{x}(f(X_{t-s})\ |\ T-s< T_0)$, we deduce that
 \begin{align*}
   R_{u,s}^T (R_{s,t}^T f)(x) &= \frac{\E_x(\E_{X_{s-u}}(f(X_{t-s})\1_{T-s < T_0}))}{\P_x(T-u<T_0)}\\
   &= \frac{\E_x(f(X_{t-s+(s-u)})\1_{T-s+(s-u) < T_0}))}{\P_x(T-u<T_0)}\\
   &=R_{u,t}^T f(x),
 \end{align*}
where we have used the Markov property a second time. Thus the
family $(R_{s,t}^T)_{0\leq s\leq t\leq T}$ is a semi-group.

\me Let us now prove that, for any $s\leq T-1$, any $x\in\N^*$ and $f\geq 0$,
\begin{equation}
\label{equation:dobrushin}
R_{s,s+1}^T f(x)\geq \frac{c_4 c_1 c_3}{2c_2} f(x_0).
\end{equation}
In other words, we prove that $\frac{c_4 c_1 c_3}{2c_2}$ is a Dobrushin coefficient,
which will allow us to show that
inequality~\eqref{equation:mixing:in-theorem:main} holds. We have
\begin{align*}
 \P_x(T-s<T_0) R_{s,s+1}^T f(x)&=\E_x\left(f(X_1) \1_{T-s<T_0}\right)\\
                             &\geq f(x_0) \P_x\left( X_1=x_0,\ T-s< T_0  \right)\\
                             &\geq f(x_0) \E_x\left(\1_{X_1=x_0}\P_{x_0}(T-s-1<T_0)\right),
\end{align*}
by the Markov property. We infer from~\eqref{equation:c4} that
$\P_{x_0}(T-s-1<T_0)\geq \frac{c_1 c_3}{2c_2} \sup_{y\in
  \N^*}\P_{y}(T-s-1<T_0)$. But 
Hypothesis~H3 yields to
\begin{equation*}
  \P_x(X_1 = x_0)\geq c_4 \P_x(1<T_0),
\end{equation*}
thus
\begin{align*}
\P_x(T-s<T_0)  R_{s,s+1}^T f(x)&\geq f(x_0) c_4 \P_x(1<T_0) \frac{c_1c_3}{2c_2}\sup_{y\in
       \N^*} \P_y(T-s-1<T_0)\\
   &\geq \frac{c_4 c_1 c_3}{2c_2} f(x_0) \P_x(T-s<T_0),
\end{align*}
by the Markov property. Finally~\eqref{equation:dobrushin} holds.

 \me We are now able to prove
 inequality~\eqref{equation:mixing:in-theorem:main}. For any
 orthogonal probability measures $\mu_1,\mu_2$ on $\N^*$ and any $f\geq 0$, we have
 by~\eqref{equation:dobrushin}
\begin{equation}
\label{eq:ZZZ}
\mu_i R_{s,s+1}^T f \geq \frac{c_4 c_1 c_3}{2c_2} f(x_0), \ \text{for }i=1,2.
\end{equation}
Thus $\mu_i R_{s,s+1}^T-\frac{c_4 c_1 c_3}{2 c_2}\delta_{x_0}$ is a
positive measure whose weight is smaller than the constant $1-\frac{c_4 c_1 c_3}{2
  c_2}$. We deduce that
 \begin{align*}
   \|\mu_1 R_{s,s+1}^T-\mu_2 R_{s,s+1}^T \|_{TV} &\leq \|(\mu_1
   R_{s,s+1}^T - \frac{c_4 c_1
     c_3}{2c_2}\delta_{x_0})\|_{TV}+\|(\mu_2 R_{s,s+1}^T - \frac{c_4 c_1
     c_3}{2 c_2}\delta_{x_0})\|_{TV}\\
     &\leq 2(1-\frac{c_4 c_1c_3}{2c_2}) = (1-\frac{c_4 c_1c_3}{2 c_2})
   \|\mu_1-\mu_2\|_{TV}.
 \end{align*}
 If $\mu_1$ and $\mu_2$ are two different but not orthogonal
 probability measures, one can apply the previous result to the
 orthogonal probability measures
 $\frac{(\mu_1-\mu_2)_+}{(\mu_1-\mu_2)_+(\N^*)}$ and
 $\frac{(\mu_1-\mu_2)_-}{(\mu_1-\mu_2)_-(\N^*)}$. Then
 \begin{multline*}
   \|\frac{(\mu_1-\mu_2)_+}{(\mu_1-\mu_2)_+(\N^*)} R_{s,s+1}^T-\frac{(\mu_1-\mu_2)_-}{(\mu_1-\mu_2)_-(\N^*)} R_{s,s+1}^T \|_{TV}\\
   \leq (1-\frac{c_4 c_1c_3}{2c_2}) \|\frac{(\mu_1-\mu_2)_+}{(\mu_1-\mu_2)_+(\N^*)}-\frac{(\mu_1-\mu_2)_-}{(\mu_1-\mu_2)_-(\N^*)}\|_{TV}.
 \end{multline*}
But
 $(\mu_1-\mu_2)_+(\N^*)=(\mu_1-\mu_2)_-(\N^*)$ since $\mu_1(\N^*)=\mu_2(\N^*)=1$, then,
 multiplying the obtained inequality by $(\mu_1-\mu_2)_+(\N^*)$, we
 deduce that
 \begin{align*}
   \|(\mu_1-\mu_2)_+ R_{s,s+1}^T-(\mu_1-\mu_2)_- R_{s,s+1}^T \|_{TV}
   &\leq (1-\frac{c_4c_1c_3}{2c_2}) \|(\mu_1-\mu_2)_+-(\mu_1-\mu_2)_-\|_{TV}.
 \end{align*}
 Since $(\mu_1-\mu_2)_+-(\mu_1-\mu_2)_-=\mu_1-\mu_2$, we obtain
 \begin{equation*}
   \|\mu_1 R_{s,s+1}^T-\mu_2 R_{s,s+1}^T \|_{TV}\leq (1-\frac{c_4c_1c_3}{2 c_2}) \|\mu_1-\mu_2\|_{TV}.
 \end{equation*}
 In particular, using the semigroup property of $(R_{s,t}^T)_{s,t}$,
 we deduce that, for any $x,y\in\N^*$,
 \begin{eqnarray*}
   \|\delta_x R_{0,T}^T - \delta_y R_{0,T}^T\|_{TV}
   &=&\|\delta_x R^T_{0,T-1} R_{T-1,T}^T - \delta_y R_{0,T-1}^T R_{T-1,T}^T\|_{TV}\\
           &\leq& \left(1-\frac{c_4c_1c_3}{2 c_2}\right)\|\delta_x R_{0,T-1}^T - \delta_y R_{0,T-1}^T\|_{TV}\\
           &\leq& 2 \left(1-\frac{c_4c_1c_3}{2c_2}\right)^{[T]},
 \end{eqnarray*}
 by induction, where $[T]$ denotes the integer part of
 $T$. Inequality~\eqref{equation:mixing:in-theorem:main} of
 Theorem~\ref{theorem:main} is thus proved for any pair of initial
 probability measures $(\delta_{x},\delta_y)$, with $(x,y)\in
 \N^*\times \N^*$.
 
 \me Let us now prove that the inequality extends to any couple of
 initial probability measures. Let $\mu$ be a probability measure on
 $\N^*$ and $x\in \N^*$. We have
\begin{align*}
  \|\P_{\mu}&(X_T\in\cdot|T<T_0)-\P_{x}(X_T\in\cdot|T<T_0)\|_{TV}\\
  &=\frac{1}{\P_{\mu}(T<T_0)}\|\P_{\mu}(X_T\in\cdot)-\P_{\mu}(T<T_0)\P_{x}(X_T\in\cdot|T<T_0)\|_{TV}\\
  &\leq
  \frac{1}{\P_{\mu}(T<T_0)} \sum_{y\in\N^*}
  \mu(y)\|\P_{y}(X_T\in\cdot)-\P_{y}(T<T_0)\P_x(X_T\in\cdot|T<T_0)\|_{TV}\\
  &\leq   \frac{1}{\P_{\mu}(T<T_0)} \sum_{y\in\N^*}
  \mu(y)\P_{y}(T<T_0)\|\P_{y}(X_T\in\cdot|T<T_0)-\P_x(X_T\in\cdot|T<T_0)\|_{TV}\\
  &\leq  \frac{1}{\P_{\mu}(T<T_0)} \sum_{y\in\N^*}
  \mu(y)\P_{y}(T<T_0) 2 (1-\frac{c_4c_1c_3}{2c_2})^{[T]}\\
  &\leq 2 (1-\frac{c_4c_1c_3}{2 c_2})^{[T]}.
\end{align*}
The same procedure, replacing $\delta_x$ by any probability measure,
leads us to inequality~\eqref{equation:mixing:in-theorem:main} in
Theorem~\ref{theorem:main}.

\bi Step 3: Let us now prove that inequality~\eqref{equation:mixing:in-theorem:main} implies the
existence and uniqueness of a QSD for $X$.

\me Let us first prove the uniqueness of the QSD. If $\rho_1$ and
$\rho_2$ are two QSDs, then we have $\P_{\rho_i}(X_t\in\cdot|t<
T_0)=\rho_i$ for $i=1,2$ and any $t\geq 0$. Thus, we deduce from
inequality~\eqref{equation:mixing:in-theorem:main} that
\begin{equation*}
  \|\rho_1-\rho_2\|_{TV}\leq 2 (1-\frac{c_4c_1c_3}{2 c_2})^{[t]},\ \forall t\geq 0,
\end{equation*}
which yields to $\rho_1=\rho_2$.

\me Let us now prove the existence of a QSD. By
\cite[Proposition~1]{Meleard2011}, this is equivalent to prove the
existence of a QLD for $X$ (see the introduction). Thus it is sufficient to
prove that there exists a point $x\in\N^*$ such that
$\P_x(X_t\in\cdot|t<T_0)$ converges when $t$ goes to infinity.

\me Let $x\in\N^*$ be any point in $\N^*$. We have, for all $s,t\geq 0$,
\begin{align*}
  \|\P_x(X_t\in\cdot|t<T_0)&-\P_x(X_{t+s}\in\cdot|t+s<T_0)\|_{TV}\\
  &=\|\P_x(X_t\in\cdot|t<T_0)-\P_{\delta_x R_{s,t+s}^{t+s}}(X_{t}\in\cdot|t<T_0)\|_{TV}\\
  &\leq 2 \left(1-\frac{c_4c_1c_3}{2c_2}\right)^{[t]}
  \xrightarrow[s,t\rightarrow+\infty]{} 0.
\end{align*}
Thus any sequence $(\P_x(X_t\in\cdot|t<T_0))_{t\geq 0}$ is a Cauchy
sequence for the total variation norm. But the space of probability
measures on $\N^*$ equipped with the total variation norm is complete,
so that $\P_x(X_t\in\cdot|t<T_0)$ converges when $t$ goes to
infinity.

\me Finally, we have proved that there exists a unique quasi-stationary
distribution $\rho$ for $X$. The last assertion of Theorem~\ref{theorem:main} is
proved as follows: for any probability measure $\mu$ on $\N^*$, we have
\begin{align*} 
  \left\| \P_{\mu}(X_t\in\cdot|t<T_0)-\rho \right\|_{TV}&=\left\|
  \P_{\mu}(X_t\in\cdot|t<T_0)-\P_{\rho}(X_t\in\cdot|t<T_0)\right\|_{TV}\\
  &\leq
  2 \left(1-\frac{c_4c_1c_3}{2 c_2}\right)^{[t]} \\ 
  &\xrightarrow[t\rightarrow
    +\infty]{} 0.
\end{align*}

\bi This concludes the proof of Theorem~\ref{theorem:main}.

\section{The birth and death process case}
\label{section:application1}

In this section, we consider \textit{birth and death processes}, which are widely used to
describe the stochastic evolution
of a population whose individuals are reproducing and dying at a rate depending on the population size. A process $X$ on $\N$ is said to be a \textit{birth and death process with absorption} if there exist two families of positive constants $(b_x)_{x\geq 1}$ and $(d_x)_{x\geq 1}$ such that the transition
rate matrix $(Q(x,y))_{x,y\in\N}$ of $X$ is given by
\begin{equation*}
  Q(x,y)=
\left\lbrace
\begin{array}{l}
b_x,\text{ if }x\geq 1\text{ and }y=x+1,\\
d_x,\text{ if }x\geq1\text{ and }y=x-1,\\
0,\text{ otherwise.}
\end{array}
\right.
\end{equation*}
 The families $(b_x)_{x\geq 1}$ and $(d_x)_{x\geq 1}$ are respectively referred to as the family of birth rates and the family of death rates. Also, one easily checks that $0$ is an absorbing point for $X$.
 
\me Applying Theorem~\ref{theorem:main}, we show that the conditional distribution of a birth and death process converges exponentially fast to a uniquely determined distribution (which is then a QSD) if and only if it admits a unique QSD. Also, as it shall be seen in the proof, Hypotheses~H1, H2 and H3 are equivalent to the uniqueness of a QSD in the birth and death case case.

\me Let us now recall that existence and uniqueness criterion for birth and death processes are well known since the works of van Doorn~\cite{vanDoorn1991} (also see Hart and Pollet~\cite{HP1996}). Indeed, setting $T_z=\inf\{t\geq 0,\ X_t=z\}$, the author proved that a birth and death process has a unique QSD if and only if
  \begin{equation*}
    S:=\sup_{x\geq 1} \E_x(T_1) < +\infty,
  \end{equation*}
where $S$ can be easily computed, since, for any $z\geq 1$,
\begin{align*}
\sup_{x\geq z} \E_x(T_z)= \sum_{k\geq z+1}\frac{1}{d_k\alpha_k}\sum_{l\geq k} \alpha_l,
\end{align*}
with
$
    \alpha_k=\left(\prod_{i=1}^{k-1} b_i\right)/\left(\prod_{i=1}^{k} d_i\right).
$
However, the spectral theory tools used to prove this result are not well suited to study the speed at which the conditional distribution converges to the quasi-stationary distribution. In particular, existing results do not provide speed of convergence to the QLD nor the set of initial distributions such that the limit~\eqref{eq:QLD} holds. As explained and described by numerical computations in~\cite{Meleard2011}, this question is of first practical importance to know whether the existence of a QSD is relevant or not for the dynamic of the process. As a consequence, the following result provide a very new insight on the quasi-limiting behaviour of birth and death processes, completing the picture offered in~\cite{vanDoorn1991}.

\begin{theorem}
\label{theorem:BD}
A birth and death process $X$ admits a unique quasi-stationary distribution if and only if there exist a constant $\gamma>0$ and a probability measure $\rho$ on $\N^*$ such that, for any initial distribution $\mu$ on $\N^*$,
  \begin{align}
  \label{eq:QLD-bis}
    \|\P_\mu(X_t\in\cdot|t<T_0)-\rho \|_{TV} \leq 2 \left(1-\gamma\right)^{[t]},\ \forall t\geq 0.
  \end{align}
  In this case, $\rho$ is the unique quasi-stationary distribution associated to $X$.
\end{theorem}

We emphasize that our proof also provides a purely probabilistic argument to the already known fact that $S<+\infty$ implies existence and uniqueness of a QSD, while earlier proofs relies on much more complex arguments based on the spectral decomposition of the rate matrix $Q$. 

\begin{proof}
Let $X$ be a birth and death process.
If~\eqref{eq:QLD-bis} holds, then $\rho$ is a QLD for $X$ starting from any initial distribution and thus it is the unique QSD for $X$.

Let us now prove that the existence and uniqueness of a quasi-stationary distribution for $X$ implies that H1, H2 and H3 hold. This will imply~\eqref{eq:QLD-bis} by Theorem 1 and thus conclude the proof of Theorem~\ref{theorem:BD}.

Since $X$ is irreducible, Hypothesis H1 is satisfied for any finite subset $K\subset \N^*$. 

Setting $x_0=1$ and $\lambda_0=b_1+d_1$, we have, for any subset $K\subset \N^*$ containing $x_0$ and for any $t\geq 0$, 
\begin{align*}
\P_{x_0}\left(X_t\in K\right)\geq \P_{x_0}\left(X_s=x_0,\,\forall s\in[0,t]\right)=e^{-\lambda_0 t}.
\end{align*}
Since the birth and death process $X$ has a unique QSD, we have $S<+\infty$ (see for instance~\cite{vanDoorn1991}). In particular, we deduce that, for any $\epsilon>0$, there exists $z_{\epsilon}\geq 1$ such that
\begin{align*}
\sup_{x\geq z_{\epsilon}} \E_x(T_{z_{\epsilon}})= \sum_{k\geq z_{\epsilon}+1}\frac{1}{\d_k\alpha_k}\sum_{l\geq k} \alpha_l \leq \epsilon.
\end{align*}
The Markov inequality thus implies that
\begin{align*}
\sup_{x\geq z_{\epsilon}}\P_x(T_{z_{\epsilon}}\geq 1)\leq \epsilon.
\end{align*}
An immediate renewal argument implies that
\begin{align*}
\sup_{x\geq z_{\epsilon}}\P_x(T_{z_{\epsilon}}\geq n)\leq \epsilon^n,\ \forall n\geq 1.
\end{align*}
As a consequence, we deduce that there exists $z_0\geq 1$ such that
\begin{align*}
\sup_{x\geq z_{0}} \E_x(e^{\lambda_0 T_{z_0}})<+\infty.
\end{align*}
In particular, setting $K=\{1,2,\ldots,z_0\}$, we deduce that
\begin{align*}
\sup_{x\in\N^*} \E_x(e^{\lambda_0 T_{K}\wedge T_0})<+\infty.
\end{align*}
Hence assumption H2 is satisfied.

Let us now prove that H3 is fulfilled. We have, for any fixed $z\in\N^*$,
\begin{align*}
\inf_{x\in \N^*} \P_x(X_1=x_0\mid T_0>1)&\geq \inf_{x\in \N^*} \P_x(X_1=x_0)\\
        &\geq \inf_{x\in \N^*} \P_x(T_{x_0}\leq 1)\P_{x_0}(X_t=1,\,\forall t\in[0,1])\\
        &\geq e^{-\lambda_0} \inf_{x\in \N^*} \P_x(T_{x_0}\leq 1)\\
        &\geq e^{-\lambda_0} \inf_{x\in \N^*} \P_x(T_z\leq 1/2)\P_z(T_{x_0}\leq 1/2),
\end{align*}
where we have used the strong, then the weak Markov property. Now, by Markov inequality, we have for any $\epsilon>0$,
\begin{align*}
\sup_{x\geq z_{\epsilon}}\P_x(T_{z_{\epsilon}}\geq 1/2)\leq 2\sup_{x\geq z_{\epsilon}}\E_x(T_{z_{\epsilon}})\leq 2\epsilon.
\end{align*} 
Choosing for instance $\epsilon=1/4$, we deduce for $z=z_{1/4}$ that
\begin{align*}
\inf_{x\in \N^*} \P_x(X_1=x_0\mid T_0>1)&\geq e^{-\lambda_0} \left(\frac{1}{2}\wedge 
\inf_{x< z} \P_x(T_z\leq 1/2)\right)\P_{z}(T_{x_0}\leq 1/2).
\end{align*}
Since $X$ is irreducible, one immediately deduces that both $\inf_{x< z} \P_x(T_z\leq 1/2)$ and $\P_{z}(T_{x_0}\leq 1/2)$ are positive. In particular, we have
\begin{align*}
\inf_{x\in \N^*} \P_x(X_1=x_0\mid T_0>1)>0,
\end{align*}
so that Hypothesis H3 is fulfilled. 

Finally, Hypotheses H1, H2 and H3 are satisfied and thus, applying Theorem~\ref{theorem:main}, there exist a constant $\gamma>0$ and a probability measure $\rho$ on $\N^*$ such that, for any initial distribution $\mu$ on $\N^*$,
  \begin{equation*}
    \|\P_\mu(X_t\in\cdot|t<T_0)-\rho \|_{TV} \leq 2 \left(1-\gamma\right)^{[t]},\ \forall t\geq 0.
  \end{equation*}
  This concludes the proof of Theorem~2.
\end{proof}

\section{A criterion on the transition rate matrix of the process}
\label{section:application2}

\me In this section, we consider a stable non-explosive process $X$  and conservative in $\N^*$, and we give a sufficient criterion on its transition rate matrix $Q$ for the
existence and uniqueness of a QSD for $X$. This result is of first theoretical importance since it applies to non-irreducible Markov processes, for which the lack of QSD-related results is patent (see for instance~\cite{vanDoorn2011}). Moreover, the assumptions of the following theorem present the interest to be accessible directly on the transition matrix expression.

\begin{theorem}
  \label{theorem:countable-state-space}
  Let $(Q(x,y))_{x,y\in\N}$ be the transition rate matrix of $X$ and assume that there exists a finite subset $K\subset \N^*$ such that
  \begin{align}
    \label{equation:theorem:countable-state-space}
    \inf_{y\in \N^*\setminus K} \left(Q(y,0) + \sum_{x\in K}  Q(y,x)\right)
    &> \sup_{y\in \N^*} Q(y,0)
  \end{align}
 and  that $P_t(x,y)>0\ \forall x,y\in
  K$.
  Then  there exists a positive constant
 $\gamma\in ]0,1[$ and $\rho \in {\cal M}_1(\N^*)$ such that, for any initial distribution $\mu$ on $\N^*$ and all $t\geq 0$,
  \begin{equation*}
    \left\|\P_{\mu}\left(X_t\in\cdot\mid t<T_0\right)-\rho \right\|_{TV}\leq 2(1-\gamma)^{[t]}.
  \end{equation*}
  In particular, $\rho$ is the unique quasi-stationary distribution associated to $X$.
\end{theorem}

\me \textit{Remark.} Our result is a generalization of~\cite[Theorem~1.1]{Ferrari2007} by Ferrari and Mari\`c, where the authors assume that
$X$ is irreducible, that
\begin{equation*}
  \overline{q}\stackrel{def}{=} \sup_{x} \sum_{y\in\N\setminus\{x\}}Q(x,y)<+\infty
\end{equation*}
and that
\begin{equation}
  \label{equation:ferrari-maric-assumption}
  \alpha\stackrel{def}{=}\sum_{x\in \N^*} \inf_{y\in \N^*\setminus x} Q(y,x) > C\stackrel{def}{=}\sup_{y\in \N^*} Q(y,0).
\end{equation}
Indeed, these assumptions imply
that there exists a finite subset $K\subset\N^*$ such that
\begin{equation*}
\inf_{y\in \N^*} \sum_{x\in K\setminus\{y\}} Q(y,x) \geq \sum_{x\in K} \inf_{y\in \N^*\setminus\{x\}} Q(y,x) > C,
\end{equation*}
and thus imply inequality~\eqref{equation:theorem:countable-state-space}.
 Moreover, we implicitly allow $\overline{q}=+\infty$ and remove the irreducibility assumption.

\begin{proof}[Proof of Theorem~\ref{theorem:countable-state-space}]

\me Let us prove that Hypotheses~H1, H2 and~H3
  hold under the assumptions of
  Theorem~\ref{theorem:countable-state-space}.

\bi By assumption, we have
$$\inf_{y\in \N^*\setminus K}\sum_{x\in K}  Q(y,x) > 0.$$
It follows that $\inf_{y\in \N^*\setminus K}\P_y(T_K\leq \frac{1}{2})>0$ and then
$$
\inf_{y\in \N^*}\P_y(T_K\leq \frac{1}{2})>0.
$$
Fix $x_0\in K$. Since $K$ is finite, we have by assumption
$$
\min_{x\in K} P_x(X_{\frac{1}{2}}=x_0)>0.
$$
Using the strong Markov property, we deduce from the two above inequalities that
\begin{equation*}
\inf_{y\in \N^*}\P_y\left(T_{x_0}\in[\frac{1}{2},1]\right)>0.
\end{equation*}
But the process is assumed to be stable, so that it remains in $x_0$
during at least a time~$\frac{1}{2}$ with positive probability. We finally deduce that
\begin{equation}
\label{eq:positivity}
\inf_{y\in \N^*}\P_y\left(X_1=x_0\right)>0,
\end{equation}
which implies Hypothesis~H3.

\bi Since $K$ is finite, for any $t\geq 0$ there exists $x^{max}_{t}\in K$ such that
$$
\P_{x^{max}_t}(t<T_0)=\max_{x\in K}\P_{x^{max}_t}(t<T_0).
$$
For $t\geq 1$, the Markov property yields to
\begin{align*}
\P_x(t<T_0) &\geq \P_x(X_1=x^{max}_t)\P_{x^{max}_{t-1}}(t-1<T_0)\\
            &\geq \P_x(X_1=x^{max}_t) \P_{x^{max}_t}(t<T_0)\\
            &\geq \min_{x',x''\in K} \P_{x'}(X_1=x'')\ \P_{x^{max}_t}(t<T_0).
\end{align*}
But $K$ is finite, thus we have by assumption
$
\min_{x',x''\in K} \P_{x'}(X_1=x'')>0.
$
And we finally deduce that
$$
\inf_{t\geq 1} \frac{\min_{x\in K} \P_x(t<T_0)}{\max_{x\in K}\P_{x}(t<T_0)}\geq \min_{x',x''\in K} \P_{x'}(X_1=x'')>0.
$$
Now, for $t\in[0,1]$, we have
\begin{align*}
\inf_{t\in[0,1]}\frac{\min_{x\in K}\P_x(t<T_0)}{\sup_{x\in K} \P_x(t<T_0)} \geq \min_{x\in K}\P_x(X_s=x,\,\forall s\in[0,1])>0.
\end{align*}
Finally, we deduce that
Hypothesis~H1 is fulfilled.

\bi Since the absorption rate of the process is uniformly bounded by $C$, we have
\begin{equation*}
  \P_{x_0}\left(X_{t-1}\in \N^*\right)\geq e^{-C(t-1)}.
\end{equation*}
By the Markov property, we deduce that
\begin{equation*}
  \P_{x_0}\left(X_t=x_0\right)\geq \inf_{y\in \N^*} \P_y\left(X_1=x_0\right) e^{-C(t-1)}.
\end{equation*}
In particular, setting $\lambda_0=C$ and using~\eqref{eq:positivity}, we deduce that the second point of
Hypothesis~H2 is fulfilled from.

\bi
   Let us now set
   \begin{equation*}
     \alpha_{K}\stackrel{def}{=}\inf_{y\in \N^*\setminus K} \left(Q(y,0) + \sum_{x\in K}  Q(y,x)\right).
   \end{equation*}
   The process jumps into $K\cup\{0\}$ from any point $x\notin
     K\cup\{0\}$ with a rate bigger than $\alpha_K$. This implies that
       $T_{K}\wedge T_0$ is uniformly bounded above by an
       exponential time of rate $\alpha_{K}$. In particular, we have
       \begin{equation*}
         \sup_{x\in \N^*}\E_x\left(e^{C T_K\wedge T_0}\right)\leq \frac{\alpha_K}{\alpha_K-C}<\infty,
       \end{equation*}
       since $\alpha_K>C$ by assumption.
       As a consequence, the first part of
       Hypothesis~H2 is also fulfilled with
       $\lambda_0=C$.

       This and Theorem~\ref{theorem:main} allows us to conclude the proof of
       Theorem~\ref{theorem:countable-state-space}.

\end{proof}

\subsection*{Acknowledgement}
This work has been partially done during the visit of D. Villemonais to
the CMM. The authors acknowledge the support from Basal-Conicyt, MATH-Amsud and the TOSCA team (INRIA Nancy -- Grand Est).

\end{document}